\documentclass[11pt]{article}
\usepackage{amsfonts}

\usepackage{graphics}
\usepackage{indentfirst}
\usepackage{cite}
\usepackage{latexsym}
\usepackage{amsmath}
\usepackage{amssymb}
\usepackage[dvips]{epsfig}
\usepackage{amscd}

\newtheorem{theorem}{Theorem}[section]

\newtheorem{lemma}[theorem]{Lemma}

\newtheorem{corollary}[theorem]{Corollary}
\newtheorem{proposition}[theorem]{Proposition}

\newcommand{\n}{\rho}

\newcommand{\bt}{\begin{theorem}}
\newcommand{\bl}{\begin{lemma}}
\newcommand{\el}{\end{lemma}}
\newcommand{\et}{\end{theorem}}

\newcommand{\la}{\label}

\newcommand{\bn}{\begin{eqnarray}}
\newcommand{\en}{\end{eqnarray}}
\newcommand{\bnn}{\begin{eqnarray*}}
\newcommand{\enn}{\end{eqnarray*}}

\newcommand{\bnnn}{\begin{eqnarray*}}
\newcommand{\ennn}{\end{eqnarray*}}
\newcommand{\ben}{\begin{enumerate}}
\newcommand{\een}{\end{enumerate}}

\newcommand{\ba}{\begin{aligned}}
\newcommand{\ea}{\end{aligned}}
\newcommand{\be}{\begin{equation}}
\newcommand{\ee}{\end{equation}}

\def\O{\mathbb{R}^N}
\def\p{\partial}
\def\norm[#1]#2{\|#2\|_{#1}}

\def\lap{\triangle}

\def\lam{\lambda}
\def\ep{\varepsilon}

\makeatletter      
\@addtoreset{equation}{section}
\makeatother       

\title{ On Formation of Singularity of Spherically Symmetric Nonbarotropic Flows}

\date{}
\author{Xiangdi H{\small UANG} }

\begin{document}
\maketitle

\begin{abstract}
We study an initial boundary value problem on a ball for
 the heat-conductive system of compressible
Navier-Stokes-Fourier equations, in particular, a criterion
of breakdown of the classical solution.
For smooth initial data away from vacuum,
it is proved that the classical solution which is spherically symmetric
 loses its regularity in a finite time if and only if the {\bf density} {\it concentrates} or {\it vanishes} or the {\bf velocity} becomes unbounded around the center. One possible situation is that a vacuum ball appears around the center and the density may concentrate on the boundary of the vacuum ball simultaneously.
\end{abstract}
\footnote[0]{2010 Mathematics Subject Classification. 35Q30, 76N10}


\section{Introduction and main results}
We are concerned with the heat-conductive system of compressible Navier-Stokes-Fourier equations which reads as
\be\la{n1}
\begin{cases} \rho_t + \nabla\cdot(\rho U) = 0,\\
 (\rho U)_t + \nabla\cdot(\rho U\otimes U) + \nabla P = \mu\lap U + (\mu + \lam)\nabla(\nabla\cdot U),\\
c_V\left((\rho\theta)_t + \nabla\cdot(\rho\theta U)\right) + P\nabla\cdot U = \kappa\lap\theta + \Psi[\nabla U],
\end{cases}
\ee
where
\be
\Psi[\nabla U] = 2\mu(\mathcal{D}(U))^2 + \lambda(\nabla\cdot U)^2,\quad \mathcal{D}(U)=\frac{\nabla U+ \nabla U^t}{2}
\ee
and $t\ge 0, x\in\Omega\subset\O(N=2,3),\ \rho=\n(t,x), U=U(t,x)$ and $\theta=\theta(t,x)$ are  the density, fluid velocity and temperature respectively.  $P=P(\rho,\theta)$ is the pressure given by a state equation
\be\la{n2}
P(\rho) = R\rho\theta .
\ee

The shear viscosity   $\mu$ , the bulk one $\lambda$ and heat conductivity $\kappa$ are constants satisfying  the physical hypothesis
\be\la{n3}
\mu,\kappa>0 ,\quad \mu+\frac{N}{2}\lam\ge 0.
\ee
The domain $\Omega$ is a bounded ball with a radius b, namely,
\be
\Omega =B_b=\{x\in\O;\ |x|\le b<\infty\}.
\ee

We study an initial boundary value problem for (1.1) with the
initial condition
\be\la{ini}
(\rho,U,\theta)(0,x) =(\rho_0,U_0,\theta_0)(x),\quad x \in \Omega,
\ee
and the boundary condition
\be\la{bdc}
U=0, \quad \frac{\p\theta}{\p n}=0\quad \ x \in \partial \Omega,\quad\mbox{$\vec{n}$ is outnormal vector.}
\ee
We are looking for the smooth spherically symmetric solution $(\rho,U)$
of the problem (1.1), (1.5),(1.6) which enjoys the form
\be
\rho(t,x)=\rho(t,|x|),\quad U(t,x)=u(t,|x|)\frac{x}{|x|},\quad\theta(t,x) = \theta(t,|x|).
\ee
Then, for the initial data to be consistent with the form (1.7), we
assume the initial data $(\rho_0,U_0)$ also takes the form
\be\la{bc-3}
\rho_0=\rho_0(|x|),\quad U_0=u_0(|x|)\frac{x}{|x|},\quad\theta_0 = \theta_0(|x|).
\ee
In this paper, we further assume the initial density is uniformly positive, that is,
\be
\rho_0=\rho_0(|x|) \ge \underline{\rho} >0, \quad x \in \Omega
\ee
for a positive constant $\underline{\rho}$.

Here, it is noted that since the assumption (1.7) implies
\be
U(t,x) + U(t,-x)=0,\quad x \in \Omega.
\ee
We necessarily have $U(t,0)=0$ (also $U_0(0)=0$) as long as classical solutions are
concerned.

There are many results about the existence of local and global strong solutions
in time
of the isentropic system of compressible Navier-Stokes equations
when the initial density is uniformly positive (refer to
\cite{Be, Itaya, Ka-1, Ka-2 , Nash, Salvi, Solo, Valli-1, Valli-2}
and their generalization\cite{Mat-1, Mat-2, Mat-3, Tani}
to the full system including the conservation law of energy). On the other hand,
for the initial density allowing vacuum,
the local well-posedness of strong solutions
of the isentropic and heat-conductive system was established by Kim\cite{Kim-1, Kim-2}.
For strong solutions with spatial symmetries,
the authors in \cite{Kim-3} proved the global existence
of radially symmetric strong solutions of the isentropic system
in an annular domain, even allowing vacuum initially.

However, it still remains open whether there exist global strong solutions
which are spherically symmetric in a ball. The main difficulties lie on
the lack of estimates of the density and velocity near the center.
In the case vacuum appears, it is worth noting that
Xin\cite{Xin} established a blow-up result which shows
that if the initial density has a compact support, then
any smooth solution to the Cauchy problem of the full system of
compressible Navier-Stokes equations without heat conduction blows up
in a finite time, see more generalization in \cite{Bu}, \cite{Xin-1} .The same blowup phenomenon occurs also for the isentropic system.
Indeed,  Zhang-Fang (\cite{ZF},Theorem 1.8) showed
that if $(\rho, U) \in C^1([0,T];H^k)\,(k > 3)$
is a spherically symmetric solution to the Cauchy problem with the compact supported initial density, then the upper limit of $T$ must be finite. To deal with large discontinuous data, Hoff \cite{Ho} established global weak solutions of the symmetric compressible heat-conductive flows. On the other hand, it's unclear whether the strong (classical) solutions
lose their regularity in a finite time
when the initial density is uniformly away from vacuum.
Therefore, it is important to study the mechanism of possible blowup
of smooth solutions, which is a main issue in this paper.

In the spherical coordinates, the original system (\ref{n1}) under the assumption (1.8) takes the form
\be\label{sym}
\left\{
\ba
& \rho_t + (\rho u)_\xi = 0,\\
& (\rho u)_t + (\rho u^2)_\xi + P_r = \nu u_{\xi r},\\
& c_V((\rho\theta)_t+(\rho u\theta)_\xi ) + Pu_\xi = \kappa\theta_{r\xi} + \nu(u_\xi)^2 -\frac{2(N-1)\mu}{r^{N-1}}(r^{N-2}u^2)_r,
\ea
\right.
\ee
where
\be
\nu = 2\mu+\lam,\quad \frac{\p}{\p_\xi} = \frac{\p}{\p r} + \frac{N-1}{r}.
\ee

Without lossing of generality, we assume $c_V=1$ and $N=3$.

Now, we consider  the following Lagrangian transformation:
\be\la{lag}
t=t,\quad h=\int_0^r\rho(t,s)\, s^{2}ds, \ \eta = (\rho r^{2})^{-1}.
\ee
Then, it follows from (1.12) that
\be
h_t = -\frac{u}{\eta},\quad r_t = u, \quad r_h = \eta,
\ee
and  the system (\ref{sym}) can be further reduced to
\be\label{sys}
\left\{
\ba
& (r^2\eta)_t = (r^{2}u)_h,\quad (\Longleftrightarrow\eta_t = u_h),\\
& u_t = r^2\left(-R\frac{\theta}{r^2\eta} + \nu\left(\frac{u_h}{r_h} + \frac{2}{r}u\right)_h\right),\\
&\left(u_t + Rr^2\left(\frac{\theta}{r^2\eta}\right)_h  = \nu r^2\left(\frac{(r^2u)_h}{r^2\eta}\right)_h\right)\\
&
\ba
\theta_t & = -R\theta\frac{(r^2u)_h}{r^2\eta} + \nu r^2\eta\left(\frac{u_h}{r_h}+\frac{2u}{r}\right)^2 - 4\mu(ru^2)_h + \kappa\left(\frac{r^2\theta_h}{r_h}\right)_h\\
& = -R\theta\frac{(r^2u)_h}{r^2\eta} + \lambda r^2\eta\left(\frac{u_h}{r_h} + \frac{2}{r}u\right)^2 + 2\mu r^2\eta\left(\frac{u_h^2}{r_h^2} + \frac{2u^2}{r^2}\right)+ \kappa\left(\frac{r^2\theta_h}{r_h}\right)_h.
\ea
\ea
\right.
\ee
The  initial boundary value problem for system (\ref{sys})
\be\la{ibvp}
\ba
& (u,\eta,\theta)(0,h) = (u_0,\eta_0,\theta_0),\quad(\eta_0>0,\theta_0>0),\\
& u(t,0) = u(t,1) = 0,\ \theta_h(t,0) = \theta_h(t,1) = 0,
\ea
\ee

where $t \ge 0$, $h \in [0, 1]$ and
\be\la{mass}
1=\int_0^b\rho_0(r)\,r^2dr=\int_0^b\rho(t,r)\,r^{2}dr,
\ee
according to the conservation of mass.
Note that
\be
r(t,0)=0,\quad r(t,1)=b.
\ee

We denote $E_0$ the initial energy
\be
E_0 =  \int_0^1\left\{\frac{u_0^2}{2} + R\left(r_0^2\eta_0-\log r_0^2\eta_0-1\right) + \left(\theta_0-\log\theta_0-1\right)\right\}dh.
\ee

Our main result is stated as follows.
\begin{theorem}\label{t1}
Assume that the initial data $(\rho_0,U_0,\theta_0)$ satisfy  {\rm (1.8),(1.9),(1.10)} and
\be
(\rho_0,U_0,\theta_0)\in H^3(\Omega).
\ee
Let $(\rho,U,\theta)$ be a classical spherically symmetric solution
to the initial boundary value problem {\rm (1.1),(1.5),(1.6),(1.7)} in $[0,T]\times\Omega$,
and $T^*$ be the upper limit of $T$, that is, the maximal time of existence of the classical solution.
Then, if $T^* < \infty$, it holds that
\be\la{ca-2}
\lim\sup_{(t,|x|)\rightarrow(T^*,0)}\left(\rho(t,|x|)+ \frac{1}{\rho}(t,|x|)+|U|(t,|x|)\right)=\infty.
\ee

\end{theorem}

\medskip

\noindent
{\bf Remark 1.1}\quad
The local existence of smooth solution with initial data as in Theorem \ref{t1}
is classical and can be found, for example, in \cite{Kim-2} and references therein.
So the maximal time $T^*$ is well defined.

\medskip

\noindent
{\bf Remark 1.2}\quad
There are several results on the blowup criterion for strong and classical solutions to the isentropic and heat-conductive system (\ref{n1}) (refer to \cite{Hxd-1, Hxd-2, Hxd-3, Hxd-full, SZ, ZF-1} and references therein).
Especially, the authors in \cite{Hxd-full} established the following Serrin-type blowup criterion:
  \be\la{bl-1}
  \limsup_{T\nearrow T^*}\left(\|\rho\|_{L^\infty(0,T;\,L^\infty)}+\|U\|_{L^r(0,T;\,L^s)}\right)=\infty,
  \ee
for any $r \in [2,\infty]$ and $s \in (3, \infty]$ satisfying
\be
\frac{2}{r}+\frac{3}{s}\le 1.
\ee

Theorem \ref{t1} asserts that the formation of singularity is only due
to the concentration or cavitation of the density and  velocity around the center.
More precisely, the density anywhere away from the center is bounded up to the maximal time. Also recall that
\be
U(t,0) = 0,\quad\mbox{for all $t<T^*$}
\ee
as far as classical solution is concerned.
It indicates the possible loss of regularity of velocity at the center.

{\bf Remark 1.3} Theorem 1.1 may be viewed as an extension of recent work in \cite{Hxd-Mat} where the authors established blowup criterion for baratropic spherically symmetric Navier-Stokes equations.
\medskip

\noindent

\section{Proof of Theorem \ref{t1}}
We only prove the case when $N=3$  since the case $N=2$ is even simpler.
Throughout of this section, we assume that
$(\rho,U,\theta)$ is a classical spherically symmetric solution with the form {\rm (1.8)}
to the initial boundary value problem {\rm (1.1),(1.5),(1.6),(1.7)} in $[0,T]\times\Omega$,
and the maximal time $T^*$, the upper limit of $T$, is finite, and we denote by $C$ generic positive  constants only depending on the initial data and the maximal time $T^*$.


By simple calculations,  we have the following estimates.
\begin{lemma}
It holds that,
\be\la{mm-1}
\int_0^1\eta dh =\int_0^1\eta_0 dh,
\ee
\be\la{mm-2}
\int_0^1r^2\eta dh = \int_0^1\left(r^2\eta\right)_0dh = \frac{b^3}{3},
\ee
\be\la{mm-3}
\int_0^1\left(\frac{1}{2}u^2 + \theta\right)dh = \int_0^1\left(\frac{1}{2}u_0^2 + \theta_0\right)dh.
\ee
\end{lemma}


Also, we have the following basic energy equality.
\begin{lemma}
It holds that
\be
\mathcal{E}(t) + \int_0^t\mathcal{V}(\tau)d\tau = \mathcal{E}(0),
\ee
\end{lemma}
where

\be\la{ent-1}
\ba
&\mathcal{E}(t) = \int_0^1\left\{\frac{u^2}{2} + R\left(r^2\eta-\log r^2\eta-1\right) + \left(\theta-\log\theta-1\right)\right\}dh,\\
&
\ba
\mathcal{V}(t) &= \int_0^1\left\{\frac{\lambda r^2\eta}{\theta}\left(\frac{u_h}{r_h} + \frac{2}{r}u\right)^2 + 2\frac{\mu r^2\eta}{\theta}\left(\frac{u_h^2}{r_h^2} + \frac{2u^2}{r^2}\right) + \kappa\frac{r^2\theta_h^2}{r_h\theta^2}\right\}dh\\
& = \int_0^1\left\{(\lambda+\frac{2}{3}\mu)\frac{ r^2\eta}{\theta}\left(\frac{u_h}{r_h} + \frac{2}{r}u\right)^2 + \frac{\mu r^2\eta}{\theta}\left(\frac{u_h}{r_h} - \frac{u}{r}\right)^2 + \kappa\frac{r^2\theta_h^2}{r_h\theta^2}\right\}dh,
\ea
\ea
\ee
or equivalently,
\be\la{ent-2}
\int_{\Omega}\rho Sdx + \int_0^t\int_{\Omega}\left(\frac{1}{\theta}\Psi[\nabla U] + \kappa\frac{|\nabla\theta|^2}{\theta^2} \right)dxds  =\int_{\Omega}\rho_0S_0dx,
\ee
where
\be
S = R\Phi(\rho^{-1}) + \Phi(\theta) + \frac{1}{2}|U|^2,\quad\Phi(s) = s -\log s - 1.
\ee

In order to prove Theorem 1.1, we can show the following stronger characterization of blowup criterion, that is,
\be\la{str-1}
 \lim\sup_{(t,h)\rightarrow(T^*-,0)}\left(\rho(t,h) + \frac{1}{\rho}(t,h) + \|U\|_{L^2(t,T^*;L^\infty(B_h))}\right)=\infty.
\ee
We argue by contradiction. For original system (1.1), assume there exists a small $r_1,\ep>0$ and a constant C, such that
\be\la{euler}
\rho(t,r) + \frac{1}{\rho}(t,r) + |U|(t,r) \le C,\quad\mbox{for $(t,r)\in (T^*-\ep,T^*)\times [0,r_1]$}.
\ee
Through the Lagrangian transformation (\ref{lag}), one immediately conclude that for system (\ref{sys}), there exists a small constant $h_1>0$, such that
\be\la{lag-1}
\rho(t,h) + \frac{1}{\rho}(t,h) + |U|(t,h) \le C,\quad\mbox{for $(t,h)\in (T^*-\ep,T^*)\times [0,h_1]$}.
\ee
Denote
\be\la{hh}
h_0 = \frac{1}{2}h_1.
\ee
Recall blowup criterion (\ref{bl-1}) by taking $r=\infty,\, s=2$, it amounts to prove
the following Proposition.

\begin{proposition}
For $h_0=\frac{1}{2}h_1$, there exists a constant C depending on $h_0$ such that
\be\la{str-1}
\rho(T,h) +\|U\|_{L^2(0,T;L^\infty(B_h^c))} \le C(h_0) \quad \mbox{for $(T,h)\in(T^*-\ep,T^*)\times [h_0,1]$}.
\ee
\end{proposition}

\medskip

 To do that, we prepare the
next lemma which gives a relationship between $r$ and $y$.
\begin{lemma}\la{le-1}
There exists a positive constant $C$ independent of T and two strict increasing function $\Theta_i:[0,1]\rightarrow [0,\infty)$  such that
  \be\la{tea-1}
  r(t,h)\ge C\Theta_1(h)
  \ee
and
\be\la{tea-2}
b^3-r(t,h)^3\ge C\Theta_2(h),
\ee
for all $0\le t<T^*$.

\end{lemma}
{\it Proof.}  For $s\ge 0$, set
\be
G(s) = s\log s - s + 1.
\ee
Obviously, $G$ is a convex function in $(0,\infty)$. By Jensen's inequality, one has
\be\la{je}
G\left(\frac{\int_{B_r}\rho dx}{|B_r|}\right)\le \frac{\int_{B_r}G(\rho)dx}{|B_r|} \Longleftrightarrow G\left(\frac{C_0h}{r^3}\right)\le \frac{C_1\int_{\Omega}G(\rho)dx}{r^3}.
\ee
Consequently, the uniform estimates for $r(t,h)$ follows immediately from Entropy inequality (\ref{ent-2}) and (\ref{je}). That is, given $h>0$ there is a strict increasing function $\Theta_1(h)$ such that
\be
r(t,h)\ge \Theta_1(h)\quad \mbox{$\Theta_1(0)=0$ and $\Theta_1(h)>0$ for $h>0$.}
\ee
By a similar step, one can conclude (\ref{tea-2}) by the following Jesen's inequality
\be
G\left(\frac{\int_{B_r^c}\rho dx}{|B_r^c|}\right)\le \frac{\int_{B_r^c}G(\rho)dx}{|B_r^c|}.
\ee

\medskip

We are now in a position to establish the pointwise estimates of the density away from the center. To do that, we first write the density in the following form. One may refer to \cite{naga} for a similar representation in one-dimensional case.
\begin{lemma}\la{le-2}
\be\la{density}
r^2\eta(t,h) = \frac{1}{\mathcal{B}(t,h)\mathcal{Y}(t,h)}\left((r^2\eta)_0(h) + \int_0^t\frac{R}{\nu}\mathcal{B}(\tau,h)\mathcal{Y}(\tau,h)\theta(\tau,h)d\tau\right).
\ee
Here
\be\la{BB}
\ba
&\mathcal{B}(t,h) = \exp\frac{1}{\nu}\bigg\{\int_{h_0}^h\frac{u_0}{r_0^2}d\xi - \int_{h_0}^h\frac{u}{r^2}d\xi - \left(\int_{h_0}^1r^2\eta dh\right)^{-1} \\
& \times \left(\int_{h_0}^1(r^2\eta)_0\int_{h_0}^h\frac{u_0}{r_0^2}d\xi dh - \int_{h_0}^1(r^2\eta)\int_{h_0}^h\frac{u}{r^2}d\xi dh + \nu\int_{h_0}^1(r^2\eta)dh -\nu\int_{h_0}^1(r^2\eta)_0dh\right) \bigg\}
\ea
\ee
and
\be\la{YY}
\ba
\mathcal{Y}(t,h) &=\exp\bigg\{\frac{1}{\nu}\int_0^t\big\{-\int_{h_0}^h\frac{2u^2}{r^3}d\xi + \left(\int_{h_0}^1r^2\eta dh\right)^{-1}\\
&\times \left\{\int_{h_0}^1\left(u^2+R\theta + r^2\eta\int_{h_0}^h\frac{2u^2}{r^3}d\xi\right) dh+ \left(\int_0^{h_0}(r^2u)_hdh\right)\left(\int_0^\tau\sigma(s,h_0)ds\right)\right\} \big\}d\tau\bigg\},
\ea
\ee
where
\be
\sigma(t,h) = -R\frac{\theta}{r^2\eta} + \nu\frac{(r^2\eta)_t}{r^2\eta}.
\ee
\end{lemma}
{\it Proof}.\
In view of (\ref{sys}), it holds
\be\la{uu-1}
\frac{1}{r^2}u_t = \left(-R\frac{\theta}{r^2\eta} + \nu\frac{(r^2\eta)_t}{r^2\eta}\right)_h\triangleq (\sigma(t,h))_h.
\ee

Thus, for $h>h_0>0$, integrating (\ref{uu-1}) over $(h_0,h)$, we deduce that
\be\la{uu-2}
\left(\int_{h_0}^h\frac{u}{r^2}d\xi\right)_t + \int_{h_0}^h\frac{2u^2}{r^3}d\xi = \sigma(t,h)-\sigma(t,h_0).
\ee
Multiplying $\frac{r^2\eta}{\nu}$ on both sides of (\ref{uu-2}) yields
\be\la{uu-3}
\frac{r^2\eta}{\nu}\left\{\left(\int_{h_0}^h\frac{u}{r^2}d\xi\right)_t + \int_{h_0}^h\frac{2u^2}{r^3}d\xi + \sigma(t,h_0)\right\} = -\frac{R}{\nu}\theta + (r^2\eta)_t,
\ee
which is
\be\la{uu-4}
(r^2\eta)_t -\frac{1}{\nu}\left\{\left(\int_{h_0}^h\frac{u}{r^2}d\xi\right)_t + \int_{h_0}^h\frac{2u^2}{r^3}d\xi + \sigma(t,h_0)\right\}(r^2\eta) = \frac{R}{\nu}\theta.
\ee
Denoted by
\be\la{uu-5}
\mathcal{A}(t,h) = -\frac{1}{\nu}\left\{\left(\int_{h_0}^h\frac{u}{r^2}d\xi\right)_t + \int_{h_0}^h\frac{2u^2}{r^3}d\xi + \sigma(t,h_0)\right\}.
\ee
In view of (\ref{uu-4}) and (\ref{uu-5}), one has
\be\la{uu-6}
(r^2\eta)(t,h) = \exp\left(-\int_0^t\mathcal{A}d\tau\right)\cdot\left\{(r^2\eta)_0(h) + \int_0^t \exp\left(\int_0^\tau\mathcal{A}ds\right)\cdot\frac{R}{\nu}\theta d\tau\right\}.
\ee
On the other hand, recall (\ref{uu-1})-(\ref{uu-2})
\be\la{uu-7}
\sigma(t,h_0) = -\left(\int_{h_0}^h\frac{u}{r^2}d\xi\right)_t -\int_{h_0}^h\frac{2u^2}{r^3}d\xi -R\frac{\theta}{r^2\eta} + \nu\frac{(r^2\eta)_t}{r^2\eta}.
\ee
Multiplying $r^2\eta$ on both sides of (\ref{uu-7}), the first term becomes
\be\la{uu-8}
\ba
-(r^2\eta)\left(\int_{h_0}^h\frac{u}{r^2}d\xi\right)_t & = -\left\{(r^2\eta)\int_{h_0}^h\frac{u}{r^2}d\xi\right\}_t + (r^2\eta)_t\int_{h_0}^h\frac{u}{r^2}d\xi\\
& = -\left\{(r^2\eta)\int_{h_0}^h\frac{u}{r^2}d\xi\right\}_t + \left\{(r^2u)\int_{h_0}^h\frac{u}{r^2}d\xi\right\}_h - (r^2u)\frac{u}{r^2}.
\ea
\ee

Integrating (\ref{uu-7}) on $[h_0,1]\times[0,t]$ and taking into account the boundary conditions (\ref{ibvp}) yields
\be\la{uu-9}
\ba
\int_0^t\int_{h_0}^1r^2\eta\sigma(\tau,h_0)dhd\tau = &-\int_0^t\int_{h_0}^1\left\{(r^2\eta)\int_{h_0}^h\frac{u}{r^2}d\xi\right\}_tdhd\tau + \int_0^t\int_{h_0}^1\left\{(r^2u)\int_{h_0}^h\frac{u}{r^2}d\xi\right\}_hdhd\tau  \\
& -\int_0^t\int_{h_0}^1\left\{u^2+R\theta + (r^2\eta)\int_{h_0}^h\frac{2u^2}{r^3}d\xi\right\}dhd\tau\\
& + \nu\int_0^t\int_{h_0}^1(r^2u)_hdhd\tau.
\ea
\ee
Therefore,
\be\la{uu-10}
\ba
\int_0^t\int_{h_0}^1r^2\eta\sigma(\tau,h_0)dhd\tau & = \int_{h_0}^1\left((r^2\eta)_0\cdot\int_{h_0}^h\frac{u_0}{r_0^2}d\xi\right) dh - \int_{h_0}^1\left((r^2\eta)\cdot\int_{h_0}^h\frac{u}{r^2}d\xi\right) dh\\
& -\int_0^t\int_{h_0}^1\left\{u^2+R\theta + (r^2\eta)\int_{h_0}^h\frac{2u^2}{r^3}d\xi\right\}dhd\tau \\
& +\nu\int_{h_0}^1(r^2\eta)dh - \nu\int_{h_0}^1(r^2\eta)_0dh.
\ea
\ee

In view of boundary condition and $(\ref{sys})_1$ ,the l,h.s of (\ref{uu-10}) can be written as
\be\la{uu-111}
\ba
& \int_0^t\left(\int_{h_0}^1r^2\eta dh\right)\left(\int_0^\tau\sigma(s,h_0)ds\right)'d\tau\\
& = \left(\int_{h_0}^1r^2\eta dh\right)\left(\int_0^t\sigma(\tau,h_0)d\tau\right) - \int_0^t\left(\int_{h_0}^1(r^2\eta)_t dh\right)\left(\int_0^\tau\sigma(s,h_0)ds\right)d\tau\\
& = \left(\int_{h_0}^1r^2\eta dh\right)\left(\int_0^t\sigma(\tau,h_0)d\tau\right) + \int_0^t\left(\int_0^{h_0}(r^2u)_hdh\right)\left(\int_0^\tau\sigma(s,h_0)ds\right)d\tau.
\ea
\ee
Collecting (\ref{uu-6})-(\ref{uu-111}), the proof of Lemma \ref{le-2} is completed.

We immediately have the following corollary
\begin{corollary}\la{cor-1}
Given $0<h_0<1$, for $h_0\le h\le 1$, there exists a constant C depending on $h_0$ such that
\be\la{by}
C^{-1}(h_0)\le \mathcal{B}(t,h), \mathcal{Y}(t,h)\le C(h_0),
\ee
and
\be\la{yy-1}
C^{-1}(h_0)\exp\{C^{-1}(h_0)(t-\tau)\} \le \frac{\mathcal{Y}(t,h)}{\mathcal{Y}(\tau,h)}\le C(h_0)\exp\{C(h_0)(t-\tau)\},\quad\mbox{$0\le\tau<t$}.
\ee

\end{corollary}
{\it Proof.} First,
\be\la{h0}
\int_{h_0}^1r^2\eta dh = \int_{h_0}^1\frac{1}{\rho}dh = \int_{r(h_0)}^br^2dr = \frac{b^3-r(h_0)^3}{3}.
\ee
In view of (\ref{tea-2}), one gets
\be
0<C\Theta_2(h_0)\le \int_{h_0}^1r^2\eta dh\le b^3.
\ee
Similarly,
\be\la{term-1}
0\le \int_{h_0}^1\frac{u^2}{r^3}d\xi =\int_{r(h_0)}^b\frac{\rho u^2r^2}{r^3}dr \le Cr(h_0)^{-3}\le C\Theta_1(h_0)^{-3}.
\ee
\be\la{term-2}
\ba
|\int_{h_0}^1\frac{u}{r^2}dh| &\le \min_{h\ge h_0}r(t,h)^{-2}\int_{h_0}^1|u|dh\le C\Theta_1(h_0)^{-2}\left(\int_0^1u^2dh\right)^{1/2}(1-h_0)^{1/2}\\
&\le C\Theta_1(h_0)^{-2}.
\ea
\ee

Also, one can verify
\be\la{yy-2}
\int_{h_0}^1r^2\eta\int_{h_0}^h\frac{2u^2}{r^3}d\xi dh
\ee
is bounded.

Observing that
\be\la{yy-3}
\ba
\int_{h_0}^1r^2\eta\int_{h_0}^h\frac{2u^2}{r^3}d\xi dh & =\int_{h_0}^1\bigg\{\left(\frac{r^3}{3}\int_{h_0}^h\frac{2u^2}{r^3}d\xi\right)_h- \frac{r^3}{3}\frac{2u^2}{r^3}\bigg\}dh\\
& =\frac{b^3}{3}\int_{h_0}^1\frac{2u^2}{r^3}d\xi -\int_{h_0}^1\frac{2u^2}{3}dh
\ea
\ee
is bounded from below and above by a constant $C(h_0)$.

To finish the proof, it suffices to bound
\be\la{sigma}
\int_0^t\sigma(\tau,h_0)d\tau.
\ee
Recall that the r.h.s of  (\ref{uu-10}) is bounded by some constant $C(h_0)$ and (\ref{uu-111}),  one has
\be\la{key-key}
\ba
|\int_0^t\sigma(\tau,h_0)d\tau| &  \le C(h_0) + C(h_0)\int_0^t\left(|\int_0^{h_0}(r^2u)_hdh|\right)\left(|\int_0^\tau\sigma(s,h_0)ds|\right)d\tau.\\
&
\ea
\ee
On the other hand, with the help of (\ref{lag-1}) yields that
\be\la{key-2}
\ba
\int_0^{h_0}(r^2u)_hdh & =\int_0^{h_0}(r^3\frac{u}{r})_hdh\\
& = \int_0^{h_0}3r\eta udh + \int_0^{h_0}r^3\left(\frac{u}{r}\right)_hdh\\
& \le C\max_{0\le h\le h_0}|r^2\eta|^{\frac{1}{2}}\left(\int_0^{h_0}\eta dh\right)^{\frac{1}{2}}\left(\int_0^{h_0} u^2dh\right)^{\frac{1}{2}} \\
& + \left(\int_0^{h_0}\frac{r^4}{\eta\theta}\left(\frac{u}{r}\right)_h^2\right)^{\frac{1}{2}}\left(\int_0^{h_0}r^2\eta\theta dh\right)^{\frac{1}{2}}\\
&\le C(h_0) + C(h_0)\left(\int_0^{h_0}\frac{r^4}{\eta\theta}\left(\frac{u}{r}\right)_h^2\right)^{\frac{1}{2}}\\
& = C(h_0) + C(h_0)\left(\int_0^{h_0}\frac{r^2\eta}{\theta}\left(\frac{u_h}{r_h}-\frac{u}{r}\right)^2\right)^{\frac{1}{2}}\\
& \le C(h_0) + C(h_0)\mathcal{V}(t)^{\frac{1}{2}},
\ea
\ee
where we used energy inequality, (\ref{mm-1}), (\ref{mm-3}) and (\ref{temp-1}) in the next Lemma \ref{lem-key}.

Consequently,   the desired bound for $\sigma(\tau,h_0)$ follows immediately from Gronwall's inequality and (\ref{key-key})-(\ref{key-2}).

Thus finishes the proof of Corollary \ref{cor-1}.

Hence, substituting (\ref{BB}) and (\ref{YY}) into (\ref{density}), we finally arrive at
\begin{lemma}\la{lem-key}
\be\la{den-1}
C^{-1}(h_0)\le \rho(t,h)\le C(h_0),\quad\mbox{$0<h_0\le h$},
\ee
\be\la{temp-1}
C^{-1}\le \int_0^1\theta dh\le C,
\ee
and
\be\la{velo-1}
\int_0^t\max_{h\in[h_0,1]}u^2(\tau, h)d\tau \le C(h_0).
\ee
\end{lemma}
{\bf Remark 2.1}
Theorem 1.1 follows immediately from Lemma \ref{lem-key}.

{\it Proof.} The right-hand side of (\ref{den-1}) is a direct consequence of the fact that $\theta>0$, (\ref{density}) and (\ref{by}). It remains to show the upper bound of $r^2\eta = \frac{1}{\rho}$.

{\it Step 1.} Multiplying $\frac{1}{\theta}$ on both sides of $(\ref{sys})_3$ yield
\be\la{temp-2}
\ba
(\log\theta)_t &= -R\left(\log(r^2\eta)\right)_t + (\lambda+\frac{2}{3}\mu)\frac{ r^2\eta}{\theta}\left(\frac{u_h}{r_h} + \frac{2}{r}u\right)^2 + \frac{\mu r^2\eta}{\theta}\left(\frac{u_h}{r_h} - \frac{u}{r}\right)^2\\
 & + \kappa\left(\frac{1}{\theta}\cdot\frac{r^2\theta_h}{r_h}\right)_h +\kappa\frac{\theta_h}{\theta^2}\cdot\frac{r^2\theta_h}{r_h}.
\ea
\ee
Integrating (\ref{temp-2}) over $[0,1]\times[0,t]$ and recall (1.4), (\ref{ibvp}) to obtain
\be\la{temp-3}
\big\{\int_0^1\log\theta dh - \int_0^1\log\theta_0dh\big\} \le -R
\big\{\int_0^1\log(r^2\eta)dh - \int_0^1\log(r^2\eta)_0dh\big\}.
\ee
Applying Jesen's inequality to (\ref{temp-3}) and (2.2) to arrive
\be\la{temp-4}
\ba
\log\int_0^1\theta dh & \ge\int_0^1\log\theta dh\\
&\ge \int_0^1\log\theta_0 dh - R\int_0^1\log(r^2\eta)dh + R\int_0^1\log(r^2\eta)_0dh\\
&\ge \int_0^1\big\{\log\theta_0 + R\log(r^2\eta)_0\big\}dh - R\log\int_0^1r^2\eta dh\\
& = C_1,
\ea
\ee
which gives the desired bound for $\int_0^1\theta dh$.

{\it Step 2.} Applying the mean value theorem to continuous function $\theta(t,h)$ to get
\be
\forall t>0, \, \exists h(t)\in [h_0,1], \quad s.t\quad \theta(t,h(t)) = \frac{\int_{h_0}^1\theta(t,h)dh}{1-h_0}.
\ee
Therefore, for $h\ge h_0$
\be\la{temp-5}
\ba
\theta(t,h)^{\frac{1}{2}} & = \theta(t,h(t))^{\frac{1}{2}} + \int_{h(t)}^h\frac{\theta_h}{2\theta(t,\xi)^{\frac{1}{2}}}d\xi\\
& \le C(1-h_0)^{-1}(\int_0^1\theta dh)^{\frac{1}{2}} +\left(\int_{h_0}^1\theta d\xi\right)^{\frac{1}{2}}\left(\int_{h_0}^1\frac{\theta_h^2}{4\theta^2}d\xi\right)^{\frac{1}{2}}\\
&\le C\big\{1+\max_{h\in[h_0,1]}r^2\eta(t,h)\int_{h_0}^1\frac{\theta_h^2}{r^2\eta\theta^2}d\xi\big\}^{\frac{1}{2}}\\
&\le C(h_0)\big\{1+\max_{h\in[h_0,1]}r^2\eta(t,h)\cdot\mathcal{V}(t)\big\}^{\frac{1}{2}}.
\ea
\ee
Consequently,
\be\la{temp-6}
\theta(t,h)\le C(h_0)\big\{1+\max_{h\in[h_0,1]}r^2\eta(t,h)\cdot\mathcal{V}(t)\big\},\quad\mbox{$h_0\le h\le 1$.}
\ee
{\it Step 3.}
Observing that
\be\la{density-1}
\ba
r^2\eta(t,h) & = \frac{1}{\mathcal{B}(t,h)\mathcal{Y}(t,h)}\left((r^2\eta)_0(h) + \int_0^t\frac{R}{\nu}\mathcal{B}(\tau,h)\mathcal{Y}(\tau,h)\theta(\tau,h)d\tau\right)\\
& \le \frac{1}{\mathcal{B}(t,h)\mathcal{Y}(t,h)}(r^2\eta)_0(h) \\
& + C\int_0^t\frac{\mathcal{B}(\tau,h)}{\mathcal{B}(t,h)}\cdot\frac{\mathcal{Y}(\tau,h)}{\mathcal{Y}(t,h)}\bigg\{1+\max_{h\in[h_0,1]}r^2\eta(t,h)\cdot\mathcal{V}(t)\bigg\}.
\ea
\ee
Hence,
\be\la{density-2}
\ba
\max_{h\in[h_0,1]}r^2\eta(t,h) & \le C + C\int_0^t\exp\{-\alpha(t-\tau)\}\bigg\{1+\max_{h\in[h_0,1]}r^2\eta(\tau,h)\cdot\mathcal{V}(\tau)\bigg\}d\tau\\
&\le C + C\int_0^t\exp\{-\alpha(t-\tau)\}\max_{h\in[h_0,1]}r^2\eta(\tau,h)\cdot\mathcal{V}(\tau)d\tau.
\ea
\ee
Write
\be
E(t) = \int_0^t\exp\{-\alpha(t-\tau)\}\max_{h\in[h_0,1]}r^2\eta(\tau,h)\cdot\mathcal{V}(\tau)d\tau.
\ee
One immediately has
\be\la{density-4}
\ba
& E_t \le \max_{h\in[h_0,1]}r^2\eta(t,h)\cdot\mathcal{V}(t)\le (C+CE)\mathcal{V}(t)-\alpha E,\\
& E_t + (\alpha-C\mathcal{V}(t))E\le C\mathcal{V}(t),\\
\ea
\ee
Applying Gronwall's inequality to (\ref{density-4}) yields
\be\la{density-5}
E\le C\exp\bigg\{-\int_0^t\alpha-C\mathcal{V}(\tau)d\tau\bigg\}\times C\int_0^t\exp\bigg\{\int_0^\tau\alpha-C\mathcal{V}(s)ds\bigg\}\mathcal{V}(\tau)d\tau.
\ee
The upper bound of $r^2\eta$ follows from (\ref{density-2}) and (\ref{density-5}).

{\it Step 4.} It suffices to establish a bound for the velocity. Indeed,
\be\la{veloc-1}
\ba
\left(\frac{u}{r}\right)^2(h) &\le \bigg\{\int_{h_0}^h|\left(\frac{u}{r}\right)_h|\bigg\}^2\le
\int_{h_0}^h\frac{r^4}{\theta r_h}\left(\frac{u}{r}\right)_h^2dh\cdot\int_{h_0}^h\frac{\theta r_h}{r^4}dh\\
& \le C(h_0)\int_{h_0}^1\frac{r^2\eta}{\theta}\left(\frac{u_h}{r_h}-\frac{u}{r}\right)^2dh.
\ea
\ee
That is
\be\la{velo-2}
\max_{h\in[h_0,1]}\left(\frac{u}{r}\right)^2(h) \le C(h_0)\mathcal{V}(t)\in L^1(0,T).
\ee
To conclude, (\ref{velo-1}) is a direct consequence of (\ref{velo-2}).

This finishes the proof of Lemma \ref{lem-key}.
 \bigskip

{\bf Acknowledgement}. \ The author would like to express his great thanks to Professor Matsumura for his helpful suggestion and discussion. This research   is partially supported  by Grant-in-Aid for JSPS Fellows 23-01320 and NSFC 11101392. The author would like to thank the anonymous referees to give several valuable suggestions, which improve our paper.

\medskip

\noindent
Xiangdi Huang\\
Academy of Mathematics and System Sciences\\
Chinese Academy of Sciences\\
Beijing 100190, P. R. China

\smallskip

\noindent
Department of Pure and Applied Mathematics\\
Graduate School of information Science and Technology \\
Osaka University\\
Toyonaka, Osaka 560-0043\\
Japan\\
e-mail: xdhuang@amss.ac.cn

\medskip


\end{document}